\documentclass[12pt,a4paper]{article}
\usepackage{amsfonts, latexsym, amsthm}

\def\seq#1#2#3{#1_{#2},\,\ldots,#1_{#3}}
\def\w{\widetilde}

\def\TT{{\underline{T}}}
\def\vv{{\underline{v}}}
\def\tt{{\underline{t}}}

\def\mm{{\underline{m}}}

\def\kk{\underline{k}}
\def\m{{\mathfrak{m}}}
\def\1{\underline{1}}

\def\P{\Bbb P}

\def\Z{\Bbb Z}

\def\C{\Bbb C}

\def\S{{\cal S}}
\def\O{{\cal O}}

\def\D{{\cal D}}

\newtheorem{theorem}{Theorem}

\newtheorem{proposition}{Proposition}

\newenvironment{definition}
{\smallskip\noindent{\bf Definition\/}:}{\smallskip\par}

\newenvironment{example}
{\smallskip\noindent{\bf Example\/}.}{\smallskip\par}
\newenvironment{example4}
{\smallskip\noindent{\bf Example 4\/}.}{\smallskip\par}
\newenvironment{examples}
{\smallskip\noindent{\bf Examples\/}.}{\smallskip\par}
\newenvironment{remark}
{\smallskip\noindent{\bf Remark\/}.}{\smallskip\par}
\newenvironment{remarks}
{\smallskip\noindent{\bf Remarks\/}.}{\smallskip\par}

\title{Alexander polynomials and Poincar\'e series of sets of ideals.
 }
\author{
A.~Campillo \and F.~Delgado
\thanks{First two authors were partially supported
by the grant MTM2007-64704. Address: University of Valladolid,
Dept. of Algebra, Geometry and Topology, 47011 Valladolid, Spain.
E-mail: campillo\symbol{'100}agt.uva.es,
fdelgado\symbol{'100}agt.uva.es} \and S.M.~Gusein-Zade \thanks{The
research was partially supported by the grants RFBR-007-00593,
INTAS-05-7805, and NWO-RFBR 047.011.2004.026. Address: Moscow
State University, Faculty of Mathematics and Mechanics, Moscow,
GSP-1, 119991, Russia. E-mail: sabir\symbol{'100}mccme.ru} }
\date{}
\begin{document}
\sloppy
\def\eps{\varepsilon}

\maketitle

\medskip

In \cite{Duke}, \cite{Edin}, \cite{Inv}, \cite{Comment}, there were
considered and, in some cases, computed Poincar\'e series of two
sorts of multi-index filtrations on the ring of germs of functions
on a complex (normal) surface singularity $(S,0)$ (in particular on
$(\C^2,0)$). A filtration from the first class was defined by a
curve (with several branches) on $(S,0)$. For $(S,0)=(\C^2,0)$, the
corresponding Poincar\'e series turned out to coincide with the
Alexander polynomial (in several variables) of the corresponding
(algebraic) link (when the number of branches is greater than 1;
otherwise it is equal to the corresponding zeta function, i.e. to
the Alexander polynomial in one variable divided by $1-t$).
Identifying all the variables, one gets the monodromy zeta function
of the curve singularity. A filtration from the second class (so
called divisorial one) was defined by a set of components of the
exceptional divisor of a modification of the surface singularity
$(S,0)$. Here we define a notion of the Poincar\'e series of a set
of ideals in the ring $\O_{(S,0)}$ of germs of functions on $(S,0)$,
describe its connection with Poincar\'e series of filtrations and
compute the Poincar\'e series of sets of ideals in some cases.

The discussed notion of the Poincar\'e series of a set of ideals
was inspired by the notion of the zeta function of an ideal given
in \cite{veys}. For our aim it is convenient to extend this notion
to a finite set of ideals $\{\seq{I}1r\}$, defining the
``Alexander polynomial" of this set. This is a mixture of the
notions introduced in \cite{veys} and \cite{sabbah}.

Let $(V,0)$ be a germ of an analytic space with an isolated
singular point at the origin and let $\seq I1r$ be ideals in the
ring $\O_{V,0}$ of
germs of functions on $(V,0)$.
Let $\pi : (X, D) \to (V,0)$ be a resolution of the
singularity of
$V$ and also of the set of the ideals $\{\seq I1r\}$. This means
that:
\begin{enumerate}
\item[1)] $X$ is a complex analytic manifold; \item[2)] $\pi$ is a
proper analytic map which is an isomorphism outside of the union
of the zero loci of the ideals $\seq I1r$; \item[3)] the
exceptional divisor $D = \pi^{-1}(0)$ is a normal crossing divisor
on $X$; \item[4)] for $i=1,\ldots,r$,  the lifting $I^*_i =
\pi^*I_i$ of the ideal $I_i$ to the space $X$  of the modification
is locally principal (and therefore $\pi$ is a principalization of
the ideal $I_i$); \item[5)] the union of the zero loci of the
liftings of the ideals $\seq  I1r$ to the space $X$ of the
modification is a normal crossing divisor on $X$.
\end{enumerate}
For $\kk = (\seq k1r)\in \Z^r_{\ge 0}$, let $S_{\kk}$ be the set of
points $x\in D$ such that, in some local coordinates $\seq z1n$ on
$X$ centred at the point $x$, the zero locus of the ideal $I^*_{i
x}$, $i=1,\ldots,r$, is the hypersurface $\{z_1=0\}$ with
multiplicity $k_i$, i.e. $I^*_{i x} = \langle z_1^{k_i} \rangle$.

\begin{definition}
{\it The Alexander polynomial of the set of ideals\/}  $\{I_i\}$
is the rational function (or a power series) in the variables
$\seq t1r$ given by the A'Campo type (\cite{AC}) formula:
$$
\Delta_{\{I_i\}}(\seq t1r) = \prod\limits_{\kk\in \Z^r_{\ge 0}\setminus
\{0\}}
(1- \tt^{\kk})^{-\chi(S_{\kk})}\; ,
$$
where $\chi(\cdot)$ is the Euler characteristic and
$\tt^{\,\kk} := t_1^{k_1}\cdot \ldots \cdot t_r^{k_r}$.
\end{definition}

\begin{remarks}
1. For $r=1$, this gives the notion of the zeta function
$\zeta_I(t)$ of an ideal $I$: \cite{veys}.
In \cite[Theorem 4.2]{veys} one forgot to write that $\pi$ should
be an isomorphism outside of the zero locus of the ideal $I$.
This is the reason why we formulate the notion for germs $(V,0)$
with isolated singularities. Another option is to demand that
the singular locus of $V$ is contained in the zero locus of the
ideal $I$.

2. The Alexander polynomial defined this way is not, in general, a
polynomial. It is really a polynomial in the case $V=\C^2$ and
$I_i=\langle f_i \rangle$, where $f_i=0$ are equations of the
irreducible components $C_i$ of a (reduced) plane curve singularity
$(C,0)\subset (\C^2,0)$: $C = \bigcup\limits_{i=1}^{r}C_i$, $r>1$.
In this case $\Delta_{\{I_i\}}(\seq t1r)$ coincides with the
classical Alexander polynomial in several variables of the algebraic
link $C \cap S^3_{\varepsilon} \subset S^3_{\varepsilon}$
($S^3_{\varepsilon}$ is the sphere of small radius $\varepsilon$
centred at the origin in $\C^2$): see \cite{EN}. We prefer to keep
the name in the general case.
\end{remarks}

Let $(V,0)$ be a germ of an analytic space and let $\P\O_{V,0}$ be
the projectivization of the ring $\O_{V,0}$ of germs of functions
on $(V,0)$. The following concept was described, e.g., in
\cite{IJM}.

Let $J^k_{V,0}=\O_{V,0}/\m^{k+1}$ be the space of k-jets of
functions on $(V; 0)$ ($\m$ is the maximal ideal in the ring
$\O_{V,0}$), let $\P J^k_{V,0}$ be its projectivization, and let
$\P^* J^k_{V,0}:=\P J^k_{V,0}\cup \{*\}$. Let
$\pi_k:\P\O_{V,0}\to\P^* J^k_{V,0}$ be the natural map. A subset
$A\subset \P\O_{V,0}$ is said to be cylindric if $A
=\pi_k^{-1}(B)$ for a constructible subset $B\subset\P
J^k_{V,0}\subset \P^* J^k_{V,0}$. The Euler characteristic
$\chi(A)$ of a cylindric set $A\subset \P\O_{V,0}$, $A
=\pi_k^{-1}(B)$, is defined as the Euler characteristic of the set
$B$. A function $\psi: \P\O_{V,0}\to G$ with values in an Abelian
group G is called cylindric if, for each $a\in G$, $a\ne 0$, the
set $\psi^{-1}(a)$ is cylindric. The integral
$\int_{\P\O_{V,0}}\psi d\chi$ of a cylindric function $\psi$ over
the space $\P\O_{V,0}$ with respect to the Euler characteristic is
the sum $\sum\limits_{a\in G\setminus\{0\}}\chi(\psi^{-1}(a))a$
(if this sum makes sense in the group $G$).

Let
\begin{equation}\label{filt}
\O_{V,0}=J(0)\supset J(1)\supset J(2)\supset\ldots
\end{equation}
be a filtration of the ring $\O_{V,0}$ by $\C^*$-invariant subsets
$J(k)$ (not necessarily vector subspaces) with cylindric
projectivizations $\P J(k)\subset \P\O_{V,0}$.

\begin{definition}
The Poincar\'e series of the filtration (\ref{filt}) is the series
$$
P(t)=\sum\limits_{k=0}^\infty \chi(\P J(k)\setminus\P J(k+1))\cdot
t^k\in\Z[[t]]\,.
$$
\end{definition}

This definition can be rewritten as $P(t)=\int_{\P\O_{V,0}}t^{v(g)}
d\chi$ where $t^\infty$ is assumed to be equal to zero.

\begin{example}
If all the elements $J(k)$ of the filtration are vector subspaces
of $\O_{V,0}$, then
$$
P(t)=\sum\limits_{k=0}^\infty \dim(J(k)\setminus J(k+1))\cdot
t^k\,.
$$
\end{example}

A filtration like (\ref{filt}) may be defined by the corresponding
($\C^*$-invariant) tautological function $v:\O_{V,0}\to \Z_{\ge
0}\cup\{+\infty\}$: $v(g)=\sup\{k: g\in J(k)\}$, $J(k)=\{g\in
\O_{V,0}: v(g)\ge k\}$. The function $v$ is cylindric in the sense
that, for each $k\in\Z_{\ge0}$, the set $\P J(k)$ is cylindric.

Let $\{v_1, \ldots, v_r\}$ be a set of $\C^*$-invariant cylindric
functions on $\O_{V,0}$ with values in $\Z_{\ge
0}\cup\{+\infty\}$. This set defines a multi-index filtration on
$\O_{V,0}$: for $\vv=(v_1, \ldots, v_r)\in\Z_{\ge 0}^r$,
\begin{equation}\label{m-filt}
J(\vv)=\{g\in \O_{V,0}: \vv(g)\ge\vv \}
\end{equation}
($\vv(g):=(v_1(g), \ldots, v_r(g))$, $\vv'\ge\vv''$ iff $v'_i\ge
v''_i$ for all $i$.). The sets $J(\vv)$ are not, in general, vector
subspaces of $\O_{V,0}$. The sets $\P J(\vv)$ are cylindric.

\begin{definition}
The Poincar\'e series of the filtration (\ref{m-filt}) is the
series
\begin{eqnarray}
P(t_1, \ldots, t_r)&=&\int_{\P\O_{V,0}}\tt^{\vv(g)} d\chi \nonumber
\\
&=& \sum\limits_{\vv\in\Z_{\ge 0}^r} \chi(\P
J(\vv)\setminus\bigcup\limits_{i=1}^r \P
J(\vv+\1_i))\cdot\tt^\vv\in\Z[[t_1, \ldots, t_r]]\,.\label{mPoi}
\end{eqnarray}
where $\1_i:=(0,\ldots, 0, 1, 0,\ldots, 0)$ ($1$ is at the $i$th
place), $\tt^\vv:=t_1^{v_1}\cdot\ldots\cdot t_r^{v_r}$\,.
\end{definition}

\begin{example}
If all $v_i$, $i=1,\ldots, r$, are order functions, that is
$v_i(g_1+g_2)\ge\min\{v_i(g_1), v_i(g_2)\}$, all the elements
$J(\vv)$ of the filtration are vector subspaces of $\O_{V,0}$, and
the definition (\ref{mPoi}) of the Poincar\'e series becomes
equivalent to the definition in \cite{Duke, IJM} given in terms of
dimensions.
\end{example}

Now let $(S,0)$ be a germ of a normal surface singularity. For a
function germ $g\in \O_{S,0}$ and a divisor $\gamma$ on $(S,0)$,
there is defined the intersection number $(\gamma\circ g)\in
\Z\cup\{\infty\}$. It can be defined, e.g., in the following way.
Let $\gamma=\sum\limits_i k_i \gamma_i$ where $\gamma_i$ are
divisors represented by irreducible curves $(C_i,0)$ on $(S,0)$.
Then $(\gamma\circ g)= \sum\limits_i k_i (\gamma_i\circ g)$ where
$(\gamma_i\circ g)$ is the order of zero of the function $g$ on
the curve $(C_i,0)$ in a uniformization parameter. (If $g_{\vert
C_i}\equiv0$, $(\gamma_i\circ g)$ is assumed to be equal to
$+\infty$.) If $\gamma$ is a Cartier divisor, $\gamma = \{f=0\}$,
$f\in \O_{S,0}$, we shall write $(f\circ g)$ instead of
$(\gamma\circ g)$.

Let $I$ be an ideal in the ring $\O_{S,0}$ and let $v_I(g):=\min\{
(g\circ f) : f\in I\}$. One can see that the function
$v_I:\O_{S,0}\to\Z_{\ge 0}\cup\{\infty\}$ is $\C^*$-invariant and
cylindric.

\begin{definition}
{\it The filtration corresponding to the ideal\/} $I$ is the
filtration
$$
\O_{S,0} = J_I(0)\supset J_I(1) \supset J_I(2) \supset \ldots
$$
defined by the function $v_I$:
$$
J_I(k) = \{g\in \O_{S,0} : v_I(g)\ge k \}\; .
$$
\end{definition}

Now let $\{\seq I1r\}$ be a set of ideals in the ring $\O_{S,0}$.

\begin{definition}
{\it The multi-index filtration corresponding to the set of
ideals\/} $\{\seq I1r\}$ is the filtration of $\O_{S,0}$, the
elements of which are
$$
J_{\{I_i\}} (\vv) = \bigcap\limits_{i=1}^r J_{I_i}(v_i)\;
$$
($\vv = (\seq v1r)\in \Z^r_{\ge 0}$).
\end{definition}

This filtration is defined by the set $\{v_{I_i}\}$ of the
corresponding functions. The spaces $J_{\{I_i\}} (\vv)$ are not, in
general, vector subspaces of $\O_{S,0}$. One can see that the
subsets $J_{\{I_i\}} (\vv)$ are, so called, "linear subspaces", i.e.
belong to the algebra generated by vector subspaces.

\begin{examples}
{\bf 1.} Let $(C,0)\subset (S,0)$ be a germ of a (reduced) curve on
the surface $S$, let $C = \bigcup\limits_{i=1}^r C_i$ be the
decomposition of the curve $C$ into irreducible components, let
$I_{C_i}$ be the ideal of the curve $C_i$. One has the multi-index
filtration corresponding to the set of ideals $\{I_{C_1}, \ldots,
I_{C_r}\}$. This is a filtration of the ring $\O_{S,0}$ by ideals.
For $S=\C^2$ or if all the components $C_i$ of the curve $C$ are
Cartier divisors on $(S,0)$ (this takes place, e.g., for any curve
on the rational double point of type $E_8$), this is the filtration
corresponding to the curve $C$ considered in \cite{Duke} and
\cite{Comment}. Otherwise this is not, generally speaking, the case.
For $S=\C^2$, the Poincar\'e series of this filtration coincides
with the Alexander polynomial (in $r$ variables) of the curve
$(C,0)$, i.e. of the corresponding link $C\cap
S_\varepsilon^{3}\subset S_\varepsilon^{3}$: \cite{Duke}.

{\bf 2.} Let $\pi : (X,D)\to (S,0)$ be a proper modification of the
space $(S,0)$ which is an isomorphism outside of the origin, with
$X$ smooth and $D = \pi^{-1}(0)$ a normal crossing divisor on $X$.
Let $D=\bigcup\limits_{\sigma\in \Gamma} E_\sigma$ be the
representation of the exceptional divisor $D$ as the union of its
irreducible components. For $\sigma\in \Gamma$, i.e. for a component
$E_\sigma$ of the exceptional divisor $D$, let $\widetilde{L}$ be a
germ of a smooth irreducible curve on $X$ intersecting $E_\sigma$
transversally at a smooth point (i.e. not at an intersection point
with other components of the exceptional divisor $D$), let
$L=\pi(\widetilde{L})$ be the corresponding curve on $(S,0)$, let
$I_{L}\subset\O_{S,0}$ be the ideal of the curve $(L,0)$, and let
$I_{E_\sigma}$ be the ideal generated by all the ideals $I_L$ of the
described type. For $r$ chosen components $\seq E1r$ of the
exceptional divisor $D$, i.e. for $\{1,\ldots, r \}\subset \Gamma$,
this way one gets a multi-index filtration (by ideals) corresponding
to a set $\{\seq E1r\}$ of components of the exceptional divisor
$D$. Again, if $S=\C^2$ or if all curves $L$ described above are
Cartier divisors on $(S,0)$, this filtration coincides with the
divisorial filtration studied in \cite{Edin} and \cite{Inv}.
Otherwise this is not the case. For $S=\C^2$, the Poincar\'e series
of this filtration coincides with the Alexander polynomial of the
set of ideals $\{I_{E_1},\ldots, I_{E_r}\}$: \cite{Edin}.

{\bf 3.} One can consider a set of valuations $v_1, \ldots,
v_{r'}, w_1, \ldots, w_{r''}$ corresponding to components $E_1,
\ldots, E_{r'}$ of a modification $\pi:(X, D)\to(S, 0)$ and to
irreducible curve germs $(C_i,0)\subset(S,0)$, $i=1, \ldots, r''$,
respectively. If $S=\C^2$, and the modification $\pi$ is an
embedded resolution of the curve
$C=\bigcup\limits_{i=1}^{r''}C_i$, one has the following formula
for the Poincar\'e series of the coresponding multi-index
filtration (by ideals). Let $m_{\sigma\delta}$ be the minus
inverse matrix of the intersection matrix $(E_\sigma\circ
E_\delta)$ of the components $E_\sigma$ on $X$, let
$\alpha_i\in\Gamma$, $i=1, \ldots, r''$, be the number of the
component $E_{\alpha_i}$ of the exceptional divisor $D$
intersecting the strict transform of the (irreducible) curve
$C_i$, and let $\mm'_{\,\sigma}:=(m_{\sigma 1}, \ldots, m_{\sigma
r'})$, $\mm''_{\,\sigma}:=(m_{\sigma \alpha_1}, \ldots, m_{\sigma
\alpha_{r''}})$. Let $\stackrel{\circ}{E_\sigma}$ be the smooth
part of the component $E_\sigma$ in the union of preimages (total
transforms) of the curves $C_i$, i.e. $E_\sigma$ minus
intersection points with other components of the exceptional
divisor $D$ and with the strict transforms of the curves $C_i$,
$i=1, \ldots, r''$. Just as in \cite{IJM} and \cite{Edin} (see
also the proof of Theorem~{\ref{theo1}}), one can show that
$$
P(t_1, \ldots, t_{r'}, T_1, \ldots,
T_{r''})=\prod\limits_{\sigma\in \Gamma}  (1 -
\tt^{\mm'_{\,\sigma}}\TT^{\mm''_{\,\sigma}})^{-\chi(
\stackrel{\circ}{E_\sigma})}
$$
(here $\tt=(t_1, \ldots, t_{r'})$ and $\TT=(T_1,\ldots, T_{r''})$
are the variables corresponding to the valuations $v_1, \ldots,
v_{r'}$ and $w_1, \ldots, w_{r''}$ respectively).
\end{examples}

Let $\{\seq I1r\}$ be a set of ideals in $\O_{\C^2,0}$. Let $\pi:
(X,D)\to (\C^2,0)$ be a resolution of the set of ideals $\{I_i\}$.
Let $D = \bigcup\limits_{\sigma\in \Gamma} E_\sigma$, where
$E_\sigma$ are irreducible components of the exceptional divisor $D$
(each $E_\sigma$ is isomorphic to the complex projective line
$\C\P^1$). For $\sigma\in \Gamma$, i.e. for a component $E_\sigma$
of the exceptional divisor $D$, let $\stackrel{\circ}{E_\sigma}$ be
the ``smooth part" of $E_\sigma$ in the union of zero loci of the
ideals $\{I^*_i\}$, i.e. $E_\sigma$ minus intersection points with
all other components of the union. Let $k_{\sigma i}$ be the
multiplicity of the component $E_\sigma$ in the zero divisor of the
ideal $I_i$ and let $\kk_{\,\sigma}:= (k_{\sigma 1}, \ldots,
k_{\sigma r})\in \Z^r_{\ge 0}$.

\begin{theorem}\label{theo1}
One has
\begin{equation}\label{eq2}
P_{\{I_i\}}(\seq t1r) = \prod\limits_{\sigma\in \Gamma}  (1 -
\tt^{\,\kk_{\,\sigma}})^{-\chi(
\stackrel{\circ}{E_\sigma})}\; .
\end{equation}
\end{theorem}

This statement generalizes those from \cite{Duke} and \cite{Edin}
for the ideals described in Examples 1 and 2 (for $S=\C^2$).
One can see that
the right-hand  side of the equation~(\ref{eq2}) is equal to the
Alexander polynomial
$\Delta_{\{I_i\}}(\seq t1r)$ of the set of ideals $\{I_i\}$
and thus in this case the Poincar\'e
series coincides with the Alexander polynomial.

\begin{proof}
The proof essentially repeats the arguments from \cite{IJM},
\cite{Edin}. One uses the representation of the Poincar\'e series
as an integral with respect to the Euler characteristic. There is
a map from the projectivization $\P\O_{\C^2,0}$ of the ring
$\O_{\C^2,0}$ onto the space of effective divisors on
$\stackrel{\circ}{D} = \bigcup\limits_{\sigma\in
\Gamma}\stackrel{\circ}{E_\sigma}$: to a function germ $f$ on
$(\C^2,0)$ one associates the intersection of the strict transform
of the curve $\{f=0\}$ with the exceptional divisor $D$. (For any
chosen $N$ this map is defined for all function germs $f$ with
$v_{I_i}(f)\le N$, $i=1, \ldots, r$, if one makes sufficiently
many additional blowing-ups at intersection points of the
components of the union of the zero loci of the ideals $I^*_i$ on
$X$. For all additional components $E_\sigma$ of the exceptional
divisor $D$, their smooth parts $\stackrel{\circ}{E_\sigma}$ are
isomorphic to the complex projective line $\C\P^1$ without two
points. Therefore their Euler characteristics are equal to zero
and they do not contribute to the right-hand side of the
equation~(\ref{eq2}).) Proposition 2 from \cite{IJM} implies that
the preimage of a point with respect to this map is a complex
affine space and thus has the Euler characteristic equal to $1$.
The Fubini formula implies that the Poincar\'e series
$P_{\{I_i\}}(\seq t1r)$ is equal to the integral with respect to
the Euler characteristic of the monomial $\tt^{\,\vv}$ over the
space of effective divisors on $\stackrel{\circ}{D}$. Here $\vv$
is an additive function on the space of effective divisors on
$\stackrel{\circ}{D}$ (with values in $\Z^r_{\ge 0}$) equal to
$\kk_{\,\sigma}$ for a point from the component
$\stackrel{\circ}{E_\sigma}$.

The space of effective divisors on
$\stackrel{\circ}{D}$ is the direct product of the spaces of
effective divisors on the components
$\stackrel{\circ}{E_\sigma}$, $\sigma\in \Gamma$. Each of the
latter ones is the disjoint union of the symmetric powers
$S^{\ell}\stackrel{\circ}{E_\sigma}$ of the component
$\stackrel{\circ}{E_\sigma}$. Therefore
$$
P_{\{I_i\}}(\seq t1r) = \prod\limits_{\sigma\in \Gamma} \left(
\sum\limits_{\ell=0}^{\infty}
\chi( S^{\ell}\stackrel{\circ}{E_\sigma})\cdot \tt^{\ell
\kk_{\,\sigma}}\right)\;
.
$$
Now equation (\ref{eq2}) follows from the formula
$$
\sum\limits_{\ell =0}^{\infty} \chi(S^{\ell} X) t^{\ell} =
(1-t)^{-\chi(X)}\;
.
$$
\end{proof}

\begin{example4}
By the definition, for ideals
$\seq I1r\subset \O_{\C^2,0}$, their Alexander polynomial
coincides with the Alexander polynomial of their integral
closures $\seq {\overline{I}}1r$. An integrally closed ideal
$I=\overline{I}$
has a representation of the form $I = \prod\limits_{\sigma}
I_{E_\sigma}^{n_\sigma}\prod\limits_{j=1}^{r''} I_{C_j}^{m_j}$
where $I_{E_\sigma}$ is the ideal corresponding to a divisor
$E_\sigma$ of a resolution of the ideal $I$ and $I_{C_j}$ is the
ideal of an irreducible curve singularity $(C_j,
0)\subset(\C^2,0)$, $m_j>0$ (see \cite{ZS}). Let $s = \# \Gamma$
be the number of components of the exceptional divisor $D$ of the
resolution. One has $n_\sigma = \sum\limits_{\delta\in \Gamma}
m_{\sigma \delta}k_\delta$, where $k_\delta$ is the multiplicity
of the ideal $I$ on the component $E_\delta$ of the exceptional
divisor $D$, $(m_{\sigma \delta})$ is the minus inverse matrix of
the intersection matrix $(E_\sigma\circ E_\delta)$ of the
components $E_\sigma$ on $X$. Therefore one has the following
equation for the order function $v_I: \O_{\C^2,0}\setminus \{0\}
\to \Z^r_{\ge 0}\cup \{\infty\}$:
$$
v_I = \sum\limits_{\sigma} n_\sigma v_\sigma
+\sum\limits_{j=1}^{r''} m_j w_j = \sum\limits_{\sigma, \delta}
m_{\sigma \delta} k_\delta v_\sigma +\sum\limits_{j=1}^{r''} m_j
w_j
$$
where $v_\sigma$ is the order function corresponding to the
component $E_\sigma$ (Example~2), $w_j$ is the order function
corresponding to the irreducible curve $(C_j, 0)$ (Example~1). One
has the following equation
$$
P_I(t) = P(t^{n_1}, \ldots , t^{n_s}, t^{m_1}, \ldots,
t^{m_{r''}})
$$
where $P(\seq t1s, \seq T1{r''})$ is the Poincar\'e series of the
valuations $(\seq v1s, \seq w1{r''})$ corresponding to the
components $E_\sigma$ of the exceptional divisor and to the curves
$(C_j,0)$ (Example~3). Moreover, one has the following statement.

\begin{proposition}
Let $(I_1, \ldots, I_r)$ be a set of ideals in $\O_{\C^2, 0}$, and
let $I_i = \prod\limits_{\sigma}
I_\sigma^{n^{i}_\sigma}\prod\limits_{j=1}^{r''} I_{C_j}^{m_j^i}$,
$i=1,\ldots,r$, where, for any $j=1, \ldots, r''$, there exists
$i$ with $m_j^i> 0$. Then one has
$$
P_{\{I_i\}}(\seq t1r) = P(\prod\limits_{i=1}^r t_i^{n_1^{i}},
\ldots , \prod\limits_{i=1}^r t_i^{n_s^{i}}, \prod\limits_{i=1}^r
t_i^{m_1^{i}}, \ldots, \prod\limits_{i=1}^r t_i^{m_{r''}^{i}})\,.
$$
\end{proposition}
\end{example4}

In a similar way one can prove versions of the main statements
from \cite{Inv} and \cite{UAC} for ideals in the ring of functions
on a rational surface singularity or on their universal abelian
covers. Let $(S,0)$ be a rational surface singularity and let
$\pi:(X,D)\to(S,0)$ be its resolution. The link $S\cap
S^3_{\varepsilon}$ of the singularity $(S,0)$ is a rational
homology sphere and its first homology group $H=H_1(S\setminus
\{0\})$ is isomorphic to the cokernel $\Z^{\Gamma}/\mbox {Im\,} j$
of the map $j : \Z^{\Gamma}\to \Z^{\Gamma}$ defined by the
intersection matrix $(E_{\sigma}\circ E_{\sigma'})$ (the order of
the group $H$ is equal to the determinant $d$ of minus the
intersection matrix $-(E_{\sigma}\circ E_{\sigma'})$). For
$\sigma\in \Gamma$, let $h_\sigma$ be the element of the group $H$
represented by the loop in the manifold $X\setminus D\simeq
S\setminus \{0\}$ going around the component $E_\sigma$ in the
positive direction. The group $H$ is generated by the elements
$h_\sigma$ for all $\sigma\in \Gamma$.

Let $p:(\widetilde{S},0)\to (S;0)$ be the universal abelian cover
of the surface singularity $S$ (see e.g. \cite{NW}, \cite{okuma},
\cite{Ne}). The group $H$ acts on $(\widetilde{S},0)$ and the
restriction $p|_{\widetilde{S}\setminus \{0\}}$ of the map $p$ to
the complement of the origin is a (usual, nonramified) covering
$\widetilde{S}\setminus\{0\}\to S \setminus\{0\}$ with the
structure group $H$. One can lift the map $p$ to a (ramified)
covering
$p':(\widetilde{X},\widetilde{D})\to (X,D)$ where $\widetilde{X}$
is a normal surface (generally speaking, not smooth) and
$\widetilde{X}\setminus \widetilde{D} \simeq
\widetilde{S}\setminus \{0\}$:
$$
\begin{array}{ccc}
(\widetilde{X}, \widetilde{\D}) &
\stackrel{\widetilde{\pi}}{\longrightarrow}
& (\widetilde{\S},0) \\
\  \downarrow {p'} &  & {p} \downarrow \ \\
(X, \D) & \stackrel{\pi}{\longrightarrow}  & (\S,0) \\
\end{array}
$$
(one can define $\w{X}$ as the normalization of the fibre
product $X\times_{\S}\widetilde{\S}$ of the varieties $X$ and
$\widetilde{\S}$ over $\S$).

Let $R(H)$ be the ring of (virtual) representations of the group
$H$. For $\sigma\in \Gamma$, i.e. for a component $E_\sigma$ of
the exceptional divisor $D$, let $\alpha_\sigma$ be the
one-dimensional representation  $H\to \C^* = {\bf GL}(1,\C)$ of
the group $H$ defined by $\alpha_\sigma(h_\delta) = \exp (- 2\pi
\sqrt{-1} m_{\sigma \delta})$ (here the minus sign reflects the
fact that the action of an element $h\in H$ on the ring
$\O_{\widetilde{\S},0}$ is defined by $(h\cdot f)(x) =
f(h^{-1}(x))$\,). Let $\{I_1, \ldots, I_r\}$ be a set of ideals in
$\O_{S,0}$ and let $\w{I}_i = p^* I_i$
be the liftings of the ideals $I_i$ to the universal abelian cover
$\w{S}$. The corresponding multi-index filtration on the ring
$\O_{\w{S},0}$ of germs of functions on the abelian cover $\w{S}$
is an $H$-invariant one. A notion of the equivariant Poincar\'e
series of such multi-index filtration was defined in \cite{mmj}.
Similar to \cite{UAC} one obtains the following result.

\begin{theorem} If the resolution $\pi:(X,D)\to(S,0)$ is a
resolution of the set of ideals $\{I_1, \ldots, I_r\}$, one has
\begin{equation}\label{eq3}
P^H (\seq t1r) = \prod\limits_{\sigma\in \Gamma} (1-
\alpha_{\sigma}\tt^{d \kk_{\,\sigma}})^{-
\chi(\stackrel{\circ}{E}_\sigma)}\; .
\end{equation}
\end{theorem}

The sum of monomials of this series with the trivial
representation in the coefficients (i.e. those for which
$\sum_{\sigma} m_{\delta \sigma} v_{\sigma}$ is an integer for
any $\delta\in \Gamma$) with the change of variables
$t_i^d \to t_i$ is the Poincar\'e series of the filtration on
$\O_{S,0}$ corresponding to the set of ideals $\{I_i\}$ (cf.
\cite{Inv}, \cite{Comment}).

\begin{remark}
Integrally closed ideals in the ring  $\O_{S,0}$ of germs of
functions on a rational surface singularity $(S,0)$ have a
description somewhat similar to that in Example~4 for $S = \C^2$.
Let $I = \overline{I}$ be an integrally closed ideal in $\O_{S,0}$
and let $\pi : (X,0)\to (S,0)$ be a resolution of it. For
$\sigma\in\Gamma$, i.e. for a component $E_\sigma$ of the
exceptional divisor $D$, let $\widetilde{L}$ be a germ of a smooth
irreducible curve on $X$ intersecting $E_\sigma$ transversally at
a smooth point, i.e. at a point of $\stackrel{\circ}{E_\sigma}$,
and let $L=\pi(\widetilde{L})$. There exists the minimal natural
number $d_\sigma$ such that $d_\sigma L $ is a Cartier divisor on
$(S,0)$: $d_\sigma L = (g_L)$ for $g_L\in \O_{S,0}$. (The number
$d_\sigma$ is the minimal natural number such that $d_\sigma
m_{\sigma \delta}$ are integers for all $\delta\in\Gamma$ and is
equal to the order of the element $h_\sigma$ in the group $H =
H_1(S\setminus\{0\})$.) Let $I'_{\sigma}\subset \O_{S,0}$ be the
ideal generated by all germs $g_L$ of the described type.

The integrally closed ideal $I$ has a unique representation of the
form $I = \prod\limits_{\sigma} I_{\sigma}^{r_\sigma}$ for non
negative rational numbers $r_\sigma$ such that
$\sum\limits_{\sigma} (E_\sigma\circ E_\delta) r_\sigma$ are
integers for all $\delta\in \Gamma$ (see \cite{Lip}). (For
integrally closed ideals $I_1$ and $I_2$ in the ring $\O_{S,0}$,
one writes that $I_1 = I_2^{1/q}$ if and only if $I_1^{q} = I_2$.
If such an ideal $I_1$ exists, it is defined in a unique way.) As
it was mentioned, the divisorial filtration on the ring $\O_{S,0}$
corresponding to the component $E_\sigma$ is not, generally
speaking, the filtration corresponding to an ideal. However, for
the corresponding functions one has $v_{I'_\sigma} = d_\sigma
v_{E_\sigma}$. Therefore the Poincar\'e series of the filtration
corresponding to the ideal $I'_\sigma$ is obtained from the one of
the divisorial filtration corresponding to the component
$E_\sigma$ by substituting the variable $t$ by $t^{d_\sigma}$.
This explains a relation between the Poincar\'e series of the
filtration in the ring $\O_{S,0}$ corresponding to a set of ideals
and that of the divisorial filtration like in Example~4.
\end{remark}


\begin{thebibliography}{12}
\bibitem[1]{AC} A'Campo N. La fonction z\^eta d'une monodromie.
Comment. Math. Helv., v.50, 233--248 (1975).
\bibitem[2]{Duke} Campillo A., Delgado F., Gusein-Zade S.M.
The Alexander polynomial of a plane curve singularity via the ring
of functions on it. Duke Math. J., v.117, no.1, 125--156 (2003).
\bibitem[3]{IJM} Campillo A., Delgado F., Gusein-Zade S.M.
The Alexander polynomial of a plane curve singularity and
integrals with respect to the Euler characteristic. Internat. J.
Math., v.14, no.1, 47--54 (2003).
\bibitem[4]{Inv} Campillo A., Delgado F., Gusein-Zade S.M.
Poincar\'e series of a rational surface singularity.
Invent. math., v.155, 41--53 (2004).
\bibitem[5]{Comment} Campillo A., Delgado F., Gusein-Zade S.M.
Poincar\'e series of curves on rational surface singularities.
Commen. Math. Helv., v.80, no.1, 95--102 (2005).
\bibitem[6]{mmj} Campillo A., Delgado F., Gusein-Zade S.M.
On Poincar\'e series of filtrations on equivariant functions of
two variables.
Moscow Math. J., v.7, no.2, 243--255 (2007).
\bibitem[7]{UAC} Campillo A., Delgado F., Gusein-Zade S.M.
Universal abelian covers of rational surface singularities and
multi-index filtrations. Funct. Anal. and Appl., v.42, no.2,
83--88 (2008).
\bibitem[8]{CDK} Campillo A., Delgado F., Kiyek K. Gorenstein
property and symmetry for one-dimensional local Cohen--Macaulay
rings. Manuscripta Mathematica, v.83, no.3--4, 405--423 (1994).
\bibitem[9]{Edin} Delgado F., Gusein-Zade S.M.
Poincar\'e series for several plane divisorial valuations. Proc.
Edinb. Math. Soc. (2), v.46, no.2, 501--509 (2003).
\bibitem[10]{EN} Eisenbud D., Neumann W. Three-Dimensional Link Theory and
Invariants of Plane Curve Singularities. Ann. of Math. Stud. 110,
Princeton Univ. Press, Princeton, 1985.
\bibitem[11]{Lip} Lipman J. Rational singularities with
applications to algebraic surfaces and unique factorization.
Inst. Hautes Etudes Sci. Publ. Math., no.36, 195--279 (1969).
\bibitem[12]{Ne} N\'emethi A. Poincar\'e series associated with
surface singularities. ArXiv: math.AG/0710.0987.
\bibitem[13]{NW} Neumann W.D., Wahl J.
Universal abelian covers of surface singularities. In: Trends in
singularities, 181--190, Trends Math., Birkh\"auser, Basel
(2002).
\bibitem[14]{okuma} Okuma T. Universal abelian covers of
rational surface singularities.
J. London Math. Soc. (2), v.70, no.2, 307--324 (2004).
\bibitem[15]{veys} van Proeyen L., Veys W. The monodromy
conjecture for zeta functions associated to ideals in dimension
two. Preprint, Leuven University, 17 pp (2007).
\bibitem[16]{sabbah} Sabbah C. Modules d'Alexander et $
D$-modules. Duke Math. J., v.60, no.3, 729--814 (1990).
\bibitem[17]{ZS} Zariski O., Samuel P. Commutative Algebra, Vol
2. van Nostrand, Princeton (1960).
\end{thebibliography}
\end{document}